\documentclass[11pt]{amsart}

\usepackage{amsmath, epsfig, color}
\usepackage{amsthm}
\usepackage{hyperref}

\numberwithin{equation}{section}

\newtheorem{prop}{Proposition}
\newtheorem{lemma}[prop]{Lemma}
\newtheorem{thm}[prop]{Theorem}

\numberwithin{prop}{section}

\theoremstyle{definition}
\newtheorem{defn}[prop]{Definition}

\newtheorem{rmk}[prop]{Remark}
\newtheorem{ques}[prop]{Question}

\newcommand{\brs}[1]{\left| #1 \right|}

\renewcommand{\gg}{\gamma}
\newcommand{\gD}{\Delta}
\newcommand{\gd}{\delta}

\newcommand{\gl}{\lambda}

\newcommand{\ga}{\alpha}
\newcommand{\gb}{\beta}
\renewcommand{\ge}{\epsilon}
\newcommand{\N}{\nabla}

\renewcommand{\bar}[1]{\overline{#1}}
\newcommand{\del}{\partial}

\newcommand{\til}[1]{\widetilde{#1}}

\newcommand{\EE}{\mathcal E}

\newcommand{\RR}{\mathcal R}
\renewcommand{\SS}{\mathcal S}
\newcommand{\LL}{\mathcal L}

\renewcommand{\underbar}[1]{\underline{#1}}

\newcommand{\tr}{\mbox{tr}_g}

\DeclareMathOperator{\Ric}{Ric}
\DeclareMathOperator{\Rm}{Rm}

\DeclareMathOperator{\sect}{sect}

\begin{document}

\title[Complete conformal metrics of negative Ricci curvature]{Existence of
Complete
conformal metrics of negative Ricci curvature on manifolds with boundary}

\author{Matthew Gursky}
\address{Institute for Advanced Study\\
         Princeton, NJ 08544}
\address{University of Notre Dame\\
         Notre Dame, IN 46556}
\email{\href{mailto:mgursky@math.ias.edu}{mgursky@ias.edu}}

\author{Jeffrey Streets}
\address{Fine Hall\\
         Princeton University\\
         Princeton, NJ 08544}
\email{\href{mailto:jstreets@math.princeton.edu}{jstreets@math.princeton.edu}}

\author{Micah Warren}
\address{Fine Hall\\
         Princeton University\\
         Princeton, NJ 08544}
\email{\href{mailto:mww@princeton.edu}{mww@princeton.edu}}

\thanks{First author supported in part by NSF grant DMS-0800084, the Charles
Simonyi Endowment and the Friends of the Institute for Advanced Study.  Second
author supported by the National Science Foundation via DMS-0703660.  Third
author supported by the National Science Foundation via DMS-0901644.}

\subjclass{53C21, 35J60}

\keywords{Poincar\'e-Einstein, conformal geometry, boundary value problems,
differential relations}

\begin{abstract} We show that on a compact Riemannian manifold with boundary
there
exists $u \in C^{\infty}(M)$ such that, $u_{|\del M}
\equiv 0$ and $u$ solves the $\sigma_k$-Ricci problem.  In the case $k = n$
the metric has negative Ricci curvature.  Furthermore, we
show the existence of a complete conformally related metric on the interior
solving the $\sigma_k$-Ricci problem.  By adopting
results of \cite{MP}, we show an
interesting relationship between the complete metrics we construct and the
existence of Poincar\'e-Einstein metrics.  Finally we give a brief discussion of
the corresponding questions in the case of positive curvature.
\end{abstract}

\date{July 23, 2009}

\maketitle

\section{Introduction}
Let $(M^n, \del M, g)$ be a smooth Riemannian manifold with boundary.  Consider
the following basic question:

\begin{ques} Let $\SS \subset \RR$ be a subset of the space of curvature tensors
on $M$.  Does there exist a conformally related metric $\hat{g} = e^{-2u} g$
such that $u_{|\del M} \equiv 0$ and the curvature tensor of $\hat{g}$ is in
$\SS$?
\end{ques}
\noindent It can also be interesting to ask the question with no boundary
restriction on $u$, or to ask for $g$ to be complete on the interior of $M$.

More specifically, we can ask: Does there exist a conformally related metric
with positive/negative scalar
curvature, Ricci tensor, Schouten tensor, sectional curvature, curvature
operator?  On closed manifolds, the answer to all of these questions is ``no''
in full generality, due to the maximum principle.  In the other direction, we
note that due to Gromov's h-principle for open diffeomorphism
invariant
differential relations on open manifolds, the answer to all of these questions
is ``yes'' if $\SS$ is open and we allow ourselves to consider \emph{all}
metrics (see \cite{EM},
\cite{Geiges}).  The restriction to a conformal class is interesting in its own
right, and can have
implications the soft methods do not yield, such as the existence of a metric
compatible with a given almost-complex structure with certain curvature
properties.  Furthermore, our methods below produce \emph{complete} metrics of
negative Ricci curvature on manifolds with boundary, a conclusion certainly not
forthcoming from the soft methods.

Our first theorem is a solution to the Dirichlet problem for the
$\sigma_k$-Ricci curvature problem on manifolds with boundary.  Before stating
the theorem we need a definition.

\begin{defn} \label{conedefs} If $A$ is a symmetric matrix $\sigma_k(A)$ is the
$k$-th elementary symmetric polynomial in the eigenvalues of $A$.  Furthermore,
let $\Gamma_k^+$ be the connected component of the set $\{\sigma_k > 0 \}$ which
contains the positive definite cone.
\end{defn}

\begin{thm} \label{dirichlet} Given $(M^n, \del M, g)$ a manifold with boundary
and $1 \leq k \leq n$, there exists a
unique
function $w_k \in C^{\infty}(M)$ such that ${w_k}_{|\del M} \equiv 0$, $-
\Ric(e^{2 w_k} g) \in \Gamma_k^+$, and
$\sigma_k \left[ - g^{-1} \Ric(e^{2 w_k} g) \right] = e^{2 k w_k}$. 
\end{thm}

In particular observe that $\Ric(e^{2 w_n} g) < 0$.  It is important to note
here that many topological obstructions exist for
curvature
conditions on manifolds with boundary.  Specifically in \cite{Ananov} a
sphere-type theorem for manifolds with positive Ricci curvature and positive
second fundamental form is shown.  Results in a similar spirit appear in
\cite{Hsiung} where the classical Bonnet-Meyers and Cartan-Hadamard theorems are
extended to manifolds with boundary.  Further topological obstructions appear in
\cite{Ichida}, \cite{Nakae}.  An interesting geometric conclusion based on
curvature and mean curvature conditions appears in \cite{HangWang}.  A common
feature of all of these results is a (usually quite strong) hypothesis on the
second fundamental form.  It is to be emphasized that our conformal factors
result in metrics with completely unknown second fundamental form, and therefore
none of the previous obstructions can apply.

By solving the Dirichlet problem with larger and larger boundary data, we can
solve the ``infinite boundary data'' Dirichlet problem to produce complete
metrics with constant $\sigma_k$ curvature on manifolds with boundary.  The case
$k = 1$ of this result appeared in \cite{AM}.  Also, the existence of a complete
metric of negative Ricci curvature with constant $\sigma_k$-Ricci curvature was
shown in \cite{Guan2} with the assumption that the given background metric
already has negative Ricci curvature.  Negative Ricci curvature of the resulting
metric is a \emph{consequence} of our theorem in the case $k = n$.

\begin{thm} \label{complete} Given $(M^n, \del M, g)$ a manifold with boundary
and $1 \leq k \leq n$, there exists a unique function $w_k \in C^{\infty}(M
\backslash \del M)$ such that $e^{2w_k} g$ is complete, $- \Ric(e^{2 w_k} g) \in
\Gamma_k^+$, and $\sigma_k \left[ -
g^{-1} \Ric (e^{2w_k} g) \right] = e^{2 k w_k}$.  Also, if $r$ denotes distance
to $\del M$, one has
\begin{align*}
\lim_{x \to \del M} w_k + \ln r - \frac{1}{2} \ln (n-1) = 0.
\end{align*}
\end{thm}

This theorem has an interesting application to understanding the existence and
moduli of Poincar\'e-Einstein metrics.  In particular, in section 6 we adopt
results of \cite{MP} to our setting and exhibit the space of conformally compact
Poincar\'e-Einstein metrics on a given manifold with boundary as an intersection
of finitely many locally closed Banach manifolds in the space of conformally
compact metrics.  For the statement of this theorem we adopt the notation most
commonly used in the study of Poincar\'e-Einstein metrics.  The relevant
terminology and the constants $\til{\gb}_{k, n}$ are defined in section 6.

\begin{thm}  \label{CCThm1} Let $(X^{n+1}, g_+)$ be a conformally compact
manifold.  Let $\Theta_k$ denote the set of 
 conformally compact metrics on $X^{n+1}$ with $\sigma_k[-g_+^{-1} \Ric] =
\tilde{\beta}_{k,n}$.  

\vskip.1in \noindent $(i)$ Given a conformally compact metric $g_{+} = \rho^{-2}
\bar{g}$, and $1 \leq 
k \leq n+1$, there is a unique conformally compact metric $h_k = e^{2w_k}
\bar{g} \in \Theta_k$.

\vskip.1in \noindent $(ii)$  Let $\EE$ denote the space of Poincar\'e-Einstein
metrics.  Then 
\begin{align*}
 \mathcal{E} = \bigcap_{k=1}^{n+1} \Theta_k,
\end{align*}
Hence $\mathcal{E}$ is a finite intersection of locally closed Banach
submanifolds, and in
particular is 
always closed in the space of conformally compact metrics on $X^{n+1}$.  
\end{thm}

In fact, the characterization of $\EE$ is much weaker, only requiring that
$\Theta_k \cap \Theta_{n+1} \neq \emptyset$ for some $k < n+1$.  This fact is
captured by a family of nonlocal conformally invariant functions we define in
section 6.  In principle these invariant functions open a path towards proving
existence of new Poincar\'e-Einstein metrics.  Specifically, on K\"ahler
manifolds one has many natural families of conformal classes, and it may be
possible to show vanishing of this invariant for carefully chosen conformal
classes.

Here is an outline of the rest of the paper.  In section 2 we recall some
basic formulas and set up the continuity method we use to prove Theorem
\ref{dirichlet}.  In sections 3 and 4 we derive the $C^1$ and $C^2$ estimates
for
the continuity method respectively, and give the proof of Theorem
\ref{dirichlet}.  In section 5 we prove Theorem \ref{complete} and in section 6
discuss the relationship of these metrics to Poincar\'e-Einstein metrics and
prove Theorem \ref{CCThm1}. 
Finally in section 7 we conclude with a brief discussion of the case of positive
curvature.

\section{Setup for Theorem \ref{dirichlet}}
We will explicitly solve the case $k = n$, and discuss the extension to the case
$k < n$ at the end of the proof.  Fix $(M^n, \del M, g)$ a compact manifold with
boundary and let
\begin{gather*}
\rho = - \Ric.
\end{gather*}
We recall that if $\hat{g} = e^{-2 u} g$ one has
\begin{gather*}
\begin{split}
\hat{\Ric} =&\ \Ric + (n-2) \N^2 u + \gD u g + (n-2) \left( du \otimes du -
\brs{du}^2 g \right).
\end{split}
\end{gather*}
It follows that
\begin{gather}
\begin{split}
\hat{\rho} =&\ \rho - (n-2) \N^2 u - \gD u g - (n-2) \left( du \otimes du -
\brs{du}^2 g \right).
\end{split}
\end{gather}
Thus if we set $w = - u$, one has
\begin{gather*}
\begin{split}
\hat{\rho} =&\ \rho + (n-2) \N^2 w + \gD w g + (n-2) \left( \brs{dw}^2 g - dw
\otimes dw \right).
\end{split}
\end{gather*}
Given a conformal factor $w$, consider the tensor
\begin{gather*}
\begin{split}
W_t(w) := (1-t) g + t \rho + (n-2) \N^2 w + \gD w g + (n-2) \left( \brs{dw}^2 g
-
dw \otimes dw \right).
\end{split}
\end{gather*}
For the remainder of this section and the next two sections we relabel our dummy
variable $w$ as $u$.  Therefore, consider the Dirichlet boundary-value
Monge-Amp\`ere equation
\begin{gather*} 
\begin{split}
F_t(u) := \det W_t(u) - e^{2n u} =&\ 0 \qquad \mbox{ ($\star_t$)}\\
u_{|\del M} \equiv&\ 0.
\end{split}
\end{gather*}
Let $\Omega = \{t \in [0,1] | \exists u \in C^{4, \ga}(M) \mbox{ solving }
(\star_t), W_t(u) \in \Gamma_n^+ \}
 \}$.  A few
observations are immediate.  First of all, equation $(\star_0)$ has the unique
solution $w \equiv 0$, thus $\Omega$ is nonempty.  Also, by construction, it is
clear that $W_0 \in \Gamma_n^+$.  By the intermediate value theorem it follows
that $W_t \in \Gamma_n^+$ for all $t$, and in particular $W_1$ will be in
$\Gamma_n^+$, as soon as the continuity method is completed.  Therefore indeed a
solution to $(\star_1)$ is the function required for the
theorem.  We can show that
the
set of times $t$ such that $(\star_t)$ is solvable is open (Lemma
\ref{openness}), therefore the crux of the matter, as always, is showing
a-priori estimates, which we will take up in the next section.  Before proving
openness of $\Omega$ we show a general maximum principle which will be of use to
us.

\begin{prop} \label{mp} \textbf{Maximum principle} Suppose that $u$ and
$v$ are smooth sub and super solutions (respectively) to equation $(\star_t)$. 
If
$u\leq v$ on $\del M,$ then $u\leq v$ on $M$.
\begin{proof} Suppose that $u>v$ somewhere.  Let $C$ be the maximum of $u-v$ on
$M$, which
is attained at some point $x_{0}$ in the interior of $M$.  Then $w=u-C$ is
a strict subsolution to $(\star_t)$, hence at the point $x_0$ we conclude
\begin{align*}
w(x_{0})  & =v(x_{0})\\
dw(x_{0})  & =dv(x_{0})\\
F_{t}(w,dw,\N^{2}w)(x_{0})  & >F_{t}(v,dv,\N^{2}v)(x_{0}).
\end{align*}
It follows immediately that at the point $x_0$ we have
\begin{align*}
\det& \left[  (1-t)g+t\rho+(n-2)\N^{2}w+\gD w g + (n-2) \left( d w \otimes d w -
\brs{d w}^2 g \right) \right]\\ 
>&\ \det\left[
(1-t)g+t\rho+(n-2)\N^{2}v+\gD v g + (n-2) \left( d v \otimes d v - \brs{d v}^2 g
\right) \right].
\end{align*}
However, note that $v\geq w$ near  $x_{0}$, which means that
\begin{align*}
\gD v(x_{0})  & \geq \gD w(x_{0}),\\
\N^{2}v(x_{0})  & \geq \N^{2}w(x_{0}).
\end{align*}
Using these inequalities and the fact that $dw(x_0) = dv(x_0)$ we conclude that
at the point $x_0$,
\begin{align*}
\mathcal W :=&\ (n-2) \N^2 w + \gD w g + (n-2) \left( d w \otimes d w - \brs{d
w}^2 g \right)\\
\leq&\ (n-2) \N^2 v + \gD v g + (n-2) \left( d v \otimes d v - \brs{d v}^2 g
\right)\\
=:&\ \mathcal V
\end{align*}
where the matrix inequality $\mathcal W \leq \mathcal V$ has the usual
interpretation that $\mathcal V - \mathcal W$ is positive semidefinite.  We
therefore conclude that
\begin{align*}
\det \left[ (1-t) g + t \rho + \mathcal V \right] =&\ \det \left[ (1-t) g + t
\rho + \mathcal W + \left(\mathcal V - \mathcal W \right) \right]\\
\geq&\ \det \left[ (1-t) g + t \rho + \mathcal W \right]
\end{align*}
which is a contradiction, and the result follows.
\end{proof}
\end{prop}

Note that this maximum principle immediately implies uniqueness of solutions to
$(\star_t)$ for all $0 \leq t \leq 1$.  Next we observe openness of $\Omega$.

\begin{lemma} \label{openness} $\Omega$ is open in $[0, 1]$.
\begin{proof} We compute the linearized operator
\begin{align*}
F_t'(u_t)(h) =&\ T_{n-1} \left( W_t \right)^{ij} \left( (n-2) \N^2 h + \gD h
g_{i j} \right. \\
&\ \left. + (n-2) \left( 2 \left<du_t, dh
\right> g - dh \otimes du_t + du_t \otimes dh \right) \right)\\
&\ - 2 n h e^{2n u_t}.
\end{align*}
where $T_{n-1} (W_t)^{ij}$ is the $(n-1)$ Newton transformation, which is
positive definite since $W_t$ is by construction.  Thus $F_t'(u_t)$ is a
strictly elliptic operator with $C^{2, \ga}$ coefficients and negative constant
term,
and is hence invertible.  The result thus follows by the implicit function
theorem.
\end{proof}
\end{lemma}

\section{Construction of Subsolutions}

In this section we derive a subsolution to the equations $(\star_t)$ which is at
the heart of our estimates.  We begin with an auxiliary geometric construction. 
Given $(M, \del M)$, let $N = \del M$
and consider the manifold $\bar{M} = M \cup \left( N \times [0, 1] \right) /
\sim$
where for $x \in N = \del M$ we have $x \times{1} \sim x$.  One should picture
an
``exterior'' collar neighborhood of $\del M$.  Using a standard partition of
unity argument one may extend the metric $g$ to a metric $\bar{g}$ defined on
$\bar{M}$ such that $\bar{g}_{|M} = g$.  Consider a point $x_0 \in \del M$.  Fix
a point $\bar{x} \in \bar{M}
\backslash M$ in the connected component of $N$ which contains $x_0$ chosen so
that $x_0$ is the closest point to $\bar{x}$ which lies
on the boundary.  Let $r$ denote geodesic distance from $\bar{x}$.  We may
arrange things so that $d(\bar{x}, \del M) > \gd$ where $\gd$ only depends on
the background metric.  See Figure \ref{fig:figure1} below.

\begin{figure}[ht]
\begin{center} \resizebox{275pt}{!}{\input{collar2.pstex_t}}
\end{center}
\caption{Exterior collar neighborhood of $\del M$} \label{fig:figure1}
\end{figure}

Fix constants $A$ and $p$ whose exact size will be determined later, and let
\begin{align*}
\underbar{u} := A \left(\frac{1}{r^p} - \frac{1}{r(x_0)^p} \right)
\end{align*}
Our goal is to show that $\underbar{u}$ is a subsolution of $(\star_t)$ for all
$t$.  First we recall the Hessian comparison theorem.

\begin{lemma} \label{hessiancomp1} \textbf{Hessian comparison theorem} Let
$(M^n, g)$ be a complete Riemannian manifold with $\sect \geq K$.  For any point
$p \in M$ the distance function $r(x) = d(x, p)$ satisfies
\begin{align*}
\N^2 r \leq \frac{1}{n-1} H_K(r) g
\end{align*}
where
\begin{gather*}
H_K(r) =  
\begin{cases}
(n-1) \sqrt{K} \cot \left( \sqrt{K} r \right) & K > 0\\
\frac{n-1}{r} & K = 0\\
(n-1) \sqrt{|K|} \coth \left( \sqrt{|K|} r \right) & K < 0
\end{cases}
\end{gather*}
\end{lemma}

\begin{lemma} \label{hessiancomp2} Let $(\bar{M}, \bar{g})$ be the metric we
constructed above, and let $r$ denote distance from a point $\bar{x} \in \bar{M}
\backslash M$ chosen so that $d(\bar{x}, \del M) > \gd > 0$ for some fixed small
constant $\gd$.  Then there exists a
constant $C$ such that 
\begin{align*}
\N^2 r(x) \leq \frac{C}{r(x)} g
\end{align*}
holds at any point where $r$ is smooth.
\begin{proof} If $K \geq 0$ the result follows immediately from Lemma
\ref{hessiancomp1}.  Assume $K \leq 0$.  The distance of $\bar{x}$ to any point
in $M$ is bounded, therefore standard estimates on the $\coth$ function give
this result
away from a controlled ball around $\bar{x}$, which we can assume is contained
in $\bar{M} \backslash M$.
\end{proof}
\end{lemma}

\begin{lemma} \label{subsolnlemma1} For $A$ and $p$ chosen large enough with
respect to constants depending only on $g$, at any
point where $r$ is smooth we have
\begin{gather*}
F_t(\underbar{u}) > 0.
\end{gather*}
\begin{proof}
We first compute $(n-2) \N^2 v + \gD v g$.  First we compute the action of this
operator on $\underbar{u}$.  Since
the action is linear we suppress the constant $A$ and reinsert it at the end.
\begin{align*}
\N^2 \underbar{u} =&\ p (p+1) r^{- p - 2} \N r \otimes \N r - p r^{- p - 1} \N^2
r\\
\gD \underbar{u} =&\ p(p+1) r^{-p-2} \brs{\N r}^2 - p r^{-p-1} \gD r\\
\end{align*}
Since $p > 0$, applying Lemma \ref{hessiancomp2} yields
\begin{align*}
- p r^{- p - 1} \N^2 r \geq - C p r^{- p - 2} g 
\end{align*}
Since the first term in the expression for $\N^2 \underbar{u}$ above is positive
we conclude
\begin{align*}
\N^2 \underbar{u} \geq - C p r^{-p-2} g
\end{align*}
for some constant $C = C(g)$.  Similarly applying Lemma \ref{hessiancomp2}
one has
\begin{align*}
- p r^{p - 1} \gD r \geq - C p r^{- p-2}
\end{align*}
for some constant $C$.  Note that $\brs{\N r} = 1$ at $x_1$ since $r$
is smooth here.  It follows that
\begin{align*}
\gD \underbar{u} \geq&\ r^{-p-2} \left( p(p+1) - C p \right)\\
\geq&\ \frac{p^2}{2} r^{-p-2}
\end{align*}
for $p$ chosen large with respect to $C$.  In sum we can conclude, reinserting
the factor $A$,
\begin{align*}
(n-2) \N^2 \underbar{u} + \gD \underbar{u} g \geq A \frac{p^2}{2} r^{-p-2} g.
\end{align*}
It is clear then that for $p$ chosen large with respect to universal constants
and then $A$ chosen large with respect to the diameter of $g$ we have 
\begin{gather*}
(n-2) \N^2 \underbar{u} + \gD \underbar{u} g \geq \frac{p^2}{4} g.
\end{gather*}
Now choose $p$ still larger depending on the ambient Ricci curvature, i.e. so
that $p \geq  4 \sqrt{ - \min_{v \in UTM}
\rho(v,v)}$.  Observing that the gradient terms in the definition of $W_t$ are
always positive, and noting that $\underbar{u} \leq 0$, we conclude the result.
\end{proof}
\end{lemma}

\section{Proof of Theorem \ref{dirichlet}}

\begin{lemma} \label{subsolnlemma2} Given $\underbar{u}$ as in Lemma
\ref{subsolnlemma1}, for all $t \in [0, 1]$, one has $\underbar{u} \leq u_t$.
\begin{proof} Fix a $t \in [0, 1]$ and suppose that $\underbar{u} > u_t$
somewhere.  We can fix a positive constant $C$ and a point $x_1 \in M$ achieving
the maximum of $\underbar{u} - u_t$, such that $\underbar{u} - C \leq u_t$ and
$(\underbar{u} - C)(x_1) = u_t(x_1)$.  It is clear by construction that this
point must be inside of $M$.  We also claim that $\underbar{u}$, and
equivalently, $r$, must be smooth
at this point $x_1$.  Indeed, if this were not the case, at $x_1$ there would be
two
geodesics $\gg_1$, $\gg_2$ which are each minimizing from $\bar{x}$ to $x_1$. 
Suppose $d(\bar{x}, x_1) = R$.  Let $\gg_1$ be given a unit speed
parametrization in $c$.  One concludes
\begin{gather} \label{subsolnlemma205}
\lim_{c \to R^-} \N r(\gg_1(c)) \cdot \gg_1' = 1.
\end{gather}
We next claim that
\begin{gather} \label{subsolnlemma210}
\lim_{c \to R^+} \N r(\gg_1(c)) \cdot \gg_1' < 1.
\end{gather}
The argument of the following paragraph is summarized in Figure
\ref{fig:figure2}.  Fix a constant $\ge > 0$ so small that $B_{\ge}(x_1)$ is
geodesically convex.  Consider the point $\til{x}_{\ge} = \gg_1(R + \ge)$. 
Construct a new curve $\til{\gg}$ from $\bar{x}$ to $\til{x}_{\ge}$ as follows:
follow the geodesic $\gg_2$ from $\bar{x}$ to $\gg_2(R - \ge)$, then connect
$\gg_2(R - \ge)$ to $\til{x}_{\ge}$ by the unique geodesic in $B_{\ge}(x_1)$
between these two points.  Recall that $\gg_1$ and $\gg_2$ are distinct
geodesics.  In particular, by uniqueness of solutions to ODE, it follows that
$\gg_1'(R) \neq \gg_2'(R)$ since $\gg_1(R) = \gg_2(R)$.  In particular, the
triangle formed by the three points $\gg_2(R - \ge)$, $\gg_1(R) = \gg_2(R) =
x_1$, and $\gg_1(R + \ge) = \til{x}_{\ge}$ is nondegenerate.  It follows from
the Toponogov comparison theorem that $d(\gg_2(R - \ge), \til{x}_{\ge})$ is
strictly less than the sum of the lengths of the other two sides of the
triangle, with the difference given in terms of a lower bound for the curvature
of $g$.  Specifically, there exists a $\gd > 0$ depending on this lower bound
and the angles of the triangle so that 
\begin{gather*}
d(\gg_2(R - \ge), \til{x_{\ge}}) \leq \left(2 - \gd \right) \ge
\end{gather*}
(In fact, since our triangle is very small, the curvature does not need to enter
into the bound.  One can forgo the Toponogov theorem and get a bound strictly in
terms of the angles of the triangle).  Using $\til{\gg}$ as a test curve for the
distance function, it follows that
\begin{align*}
d(\bar{x}, \til{x}_{\ge}) \leq R - \ge + \left(2 - \gd\right) \ge =&\ R + \ge -
\gd \ge.
\end{align*}
Taking the limit as $\ge \to 0$, we immediately conclude that
\begin{align*}
\lim_{c \to R^+} \N r(\gg_1(c)) \cdot \gg_1' =&\ \lim_{\ge \to 0}
\frac{r(\gg_1(R+\ge)) - r(\gg_1(R))}{\ge}\\
\leq&\ \lim_{\ge \to 0} \frac{R + \ge - \gd \ge - R}{\ge}\\
<&\ 1.
\end{align*}

\begin{figure}[ht]
\begin{center} \resizebox{275pt}{!}{\input{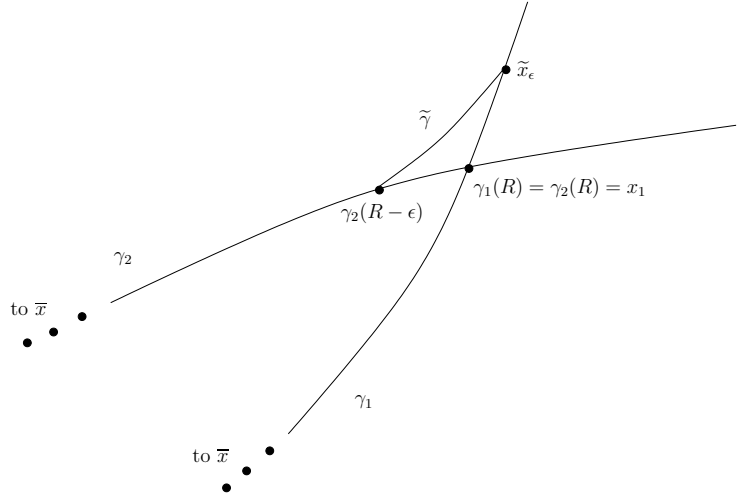}}
\end{center}
\caption{Geodesics at the cut locus} \label{fig:figure2}
\end{figure}

We now finish the argument that $\underbar{u}$ is smooth at $x_1$.  Indeed, it
follows from (\ref{subsolnlemma205}) and (\ref{subsolnlemma210}) by direct
calculation that the derivative of the function $f(c) := \underbar{u}(\gg_1(c))$
jumps
a certain positive amount at $c = R$. Considering next the smooth function
$\psi(c) := u_t(\gg_1(c))$, by assumption we have that $(\psi - f)(c)$ has a
local minimum at $c = R$.  Thus 
\begin{gather*}
\lim_{c \to R^-} f' - \psi' \leq \lim_{c \to R^+} f' - \psi'.
\end{gather*}
Since $\psi$ is smooth, we therefore conclude
\begin{align*}
\lim_{c \to R^-} f' \geq \lim_{c \to R^+} f'.
\end{align*}
This contradicts (\ref{subsolnlemma205}) and (\ref{subsolnlemma210}) since
$\frac{d\underbar{u}}{d r} < 0$.

Given that $\underbar{u}$ is smooth at $x_1$, using Lemma \ref{subsolnlemma1}
the argument of Proposition \ref{mp} applies at this point to yield the required
contradiction to the assumption that $\underbar{u} > u_t$ somewhere.  The lemma
follows.
\end{proof}
\end{lemma}

\begin{lemma} \label{subsolnlemma4} The inequality $u_t \leq 0$ holds for all $0
\leq t \leq 1$.
\begin{proof} To get this estimate we exhibit $u_t$ as a subsolution of
$(\star_0)$.  We may assume without loss of
generality that by scaling $g$ in space we have $g \geq \rho$.  It follows that
\begin{align*}
e^{2 n u_t} =&\ F_t(u_t) + e^{2n u_t}\\
=&\ \det \left( (1-t) g + t \rho + (n-2) \N^2 u_t + \dots \right)\\
\leq&\ \det \left( g + (n-2) \N^2 u_t + \dots \right)\\
=&\ F_0(u_t) + e^{2n u_t}.
\end{align*}
Therefore $u_t$ is a subsolution of the equation $F_0(u) = 0$, and the result
follows by Proposition \ref{mp}.
\end{proof}
\end{lemma}

\begin{prop} \label{NRboundaryC1} There exists a constant $C$ such that for all
$x_0 \in \del M$ and
for all $0 \leq t \leq 1$ we have $\brs{\frac{\del}{\del \nu} u_t} \leq C$ where
$\nu$ is the unit normal to $\del M$ at $x_0$.
\begin{proof} This follows immediately from Lemmas \ref{subsolnlemma2} and
\ref{subsolnlemma4} since for instance
\begin{align*}
\frac{\underbar{u}(x) - \underbar{u}(x_0)}{d(x,x_0)} \leq \frac{u(x) -
u(x_0)}{d(x,x_0)}.
\end{align*}
Our construction of $\underbar{u}$ is specific to each $x_0 \in \del M$, but it
is clear that the choice of $p$ etc. are all universally controlled, and so the
proposition follows.
\end{proof}
\end{prop}

\begin{prop} \label{C1Estimate} There exists a constant $C$ such that for all $0
\leq t \leq 1$ we
have
\begin{align*}
\brs{u_t}_{C^1} \leq C
\end{align*}
\begin{proof} We have already shown the global $C^0$ estimate and the boundary
$C^1$ estimate.  Suppose that the maximum of $\brs{\N u_t}$ occurs at a point in
the interior.  One may follow the calculation of \cite{GV} Proposition 4.1,
which is
justified at any interior point of $M$, to yield the a-priori $C^2$ estimate. 
The result follows. 
\end{proof}
\end{prop}

We now proceed with the $C^2$ estimates.  Fix $x_0 \in \del M$ and let $u$ be a
solution to $(\star_t)$ for some $0 \leq t
\leq 1$.  Suppose further that $\brs{\N_n u}  < C$.  We will use the indices
$e_{i},e_{j}$ to refer to tangent directions to $\del M$, and $e_n$ to refer
to the
unit inward normal at $x_0$.  We require separate proofs for the different types
of boundary second derivatives $\N_i \N_i u$, $\N_i \N_n u$, and $\N_n \N_n u$.

\begin{lemma} \label{TTC2Estimate} There exists a constant $C$ depending on
$\sup_{0 \leq t \leq 1} \brs{u_t}_{C^1}$ such that for all $x_0 \in \del M$, for
all $0 \leq t \leq 1$ we have
\begin{align*}
\brs{\N_i \N_j u(x_0)} < C.
\end{align*}
\begin{proof} We note first, using that $u_{| \del M} \equiv 0$,
\begin{align*}
\N_i \N_j u(x_0) =&\ - \N_n u (x_0) A(e_i, e_j)
\end{align*}
where $A$ is the second fundamental form of $\del M$.  Since $\brs{\N_n u (x_0)}
< C$ we immediately conclude the result.
\end{proof}
\end{lemma}

Next we need to bound the derivatives of the form $\N_{e_i} \N_n u$ at the
boundary.  For our given $t \in [0, 1]$, let $\LL$ denote the linearization of
$F_t$ at $u$.  As in Lemma \ref{openness} we have
\begin{gather} \label{Ldef}
\begin{split}
\LL(\psi) =&\ T_{n-1}(W_t)^{ij} \left( (n-2) \N_i \N_j \psi + \gD \psi g_{i j} +
\right.\\
&\ \qquad \left. (n-2) \left( 2 \left< \N \psi, \N u \right> g_{i j} - \N_i \psi
\otimes \N_j u - \N_i u \otimes \N_j \psi \right) \right)\\
&\ - 2 n \psi e^{2n u}.
\end{split}
\end{gather}
Fix a point $x_0 \in \del M$ and let $B_{\gd}$ be the ball of some small radius
$\gd > 0$ around $x_0$.  Pick coordinates in $B_{\gd}$ so that $\del M$ is the
plane $x_n = 0$, and let $\{e_i, e_n \}$ be the corresponding coordinate vector
fields.  Fix
some
$\ga$ and consider the function $\phi = e_{\ga} u_t$ defined in $B_{\gd}$.  Note
that $\phi_{|\del M} = 0$.  We
aim to apply a maximum principle argument similar to the $C^1$ boundary estimate
to yield a bound on the normal derivative of $\phi$, which will yield the
required bound.  The first step is to bound the action of $\LL$ on $\phi$.

\begin{lemma} \label{TNlemma1} Using the notation above, there exists a constant
$C$ such that
\begin{align*}
\brs{\LL(\phi)} \leq C \left(1 + \sup_{0 \leq t \leq 1} \brs{u_t}_{C^1} \right)
\sum F^{ii}.
\end{align*}
\begin{proof} Differentiating equation $(\star_t)$ with respect to $e_{\ga}$
yields
\begin{align*}
0 =&\ F^{ij} \left( t \N_{\ga} \rho_{ij} + (n-2) \N_{\ga} \N_i \N_j u + \N_{\ga}
\gD u g_{i j} \right.\\
&\ \qquad \left. + (n-2) \left( 2 \left< \N_{\ga} \N u, \N u \right> g - \N_\ga
\N_i u \otimes \N_j u - \N_i u \otimes \N_{\ga} \N_j u \right) \right)\\
&\ - 2 n \N_{\ga} u e^{2n u}.
\end{align*}
Commuting derivatives we conclude
\begin{align*}
\N_{\ga} \N_i \N_j u =&\ \N_i \N_{\ga} \N_j u + \Rm * \N u\\
=&\ \N_i \N_j \N_{\ga} u + \Rm * \N u\\
=&\ \N_i \N_j \phi + \Rm * \N u.
\end{align*}
Similarly we have
\begin{align*}
\N_{\ga} \gD u =&\ \gD \phi + \Rm * \N u.
\end{align*}
Combining these calculations yields
\begin{align*}
\LL(\phi) = F^{ij} \left( t \N_{\ga} \rho_{ij} + \left(\Rm * \N u \right)_{ij}
\right).
\end{align*}
Applying the $C^1$ bound we immediately conclude
\begin{align*}
\brs{\LL(\phi)} \leq C \brs{F^{ij}} \left(1 + \brs{\N u} \right) \leq C
\brs{F^{ij}}
\end{align*}
where $\brs{F^{ij}}$ refers to the matrix norm.  Note that $F^{ij}$ is positive
definite, and therefore its norm is dominated by a dimensional constant times
its trace.  The result follows.
\end{proof}
\end{lemma}

\begin{lemma} \label{TNC2Estimate} There exists a constant $C$ depending on
$\sup_{0 \leq t \leq 1} \brs{u_t}_{C^1}$ such that for all $x_0 \in \del M$, for
all $0 \leq t \leq 1$ we have
\begin{align*}
\brs{\N_{i} \N_{n} u_t(x_0)} \leq C.
\end{align*}
\begin{proof} Fix a point $x_0 \in \del M$.  Choose a small constant $\gd <
\frac{1}{2}$ so that $B_{\gd}(x_0)$ is a geodesically convex ball.  Furthermore
choose $\bar{x}$ such that $\bar{x} \in B_{\frac{\gd}{4}}(x_0)$.  Following the
construction in section 3, let
\begin{align*}
\underbar{u} = \frac{1}{r^p} - \frac{1}{r(x_0)^p}.
\end{align*}
Note that by construction we have that $\underbar{u}$ is smooth in $U := M \cap
B_{\gd}(x_0)$.  We want to compute $\LL(\underbar{u})$.  Following the
calculations
of section 3 we conclude
\begin{align*}
(n-2) \N^2 \underbar{u} + \gD \underbar{u} g \geq \frac{p^2}{4} r^{-p-2} g.
\end{align*}
Consider one of the linear terms in $\LL$ acting on $\underbar{u}$.  In
particular we note that
\begin{align*}
\left< \N \underbar{u}, \N u \right> g \leq&\ p r^{-p-1} \brs{\N u} g\\
\leq&\  C p r^{-p-2} g
\end{align*}
using that $r \leq 1$ and $\brs{\N u} \leq C$.  All of the linear terms are
bounded thusly.  Furthermore, $\underbar{u} \leq 0$, so we can throw away the
constant term and conclude
\begin{align*}
\LL(\underbar{u}) \geq \frac{p^2}{8} r^{-p-2} \sum F^{ii}
\end{align*}
when $p$ is chosen large enough with respect to fixed constants.  It follows
that if we choose $p$ larger still and note that $r \leq 1$, Lemma
\ref{TNlemma1} yields
\begin{align*}
\LL( \phi - \underbar{u} ) < \left( C - \frac{p^2}{8} \right) \sum F^{ii} < 0.
\end{align*}
Thus by the maximum principle we conclude that the minimum of $\phi -
\underbar{u}$ occurs on the boundary of $U$.  It remains to check these boundary
values.  There are two components of $U$ to check.  First we have $U \cap \del
M$.  Here $\phi \equiv 0$ and $\underbar{u} \leq 0$ with $\underbar{u} = 0$ at
$x_0$.  Next we consider the component $U \cap \del B_{\gd}(x_0)$.  Here using
the $C^1$ estimate we have $\phi \geq - C$ for some controlled constant $C$. 
Since $\bar{x} \in B_{\frac{\gd}{4}}(x_0)$, it follows that for $x \in U \cap
\del B_{\gd}(x_0)$ one has
\begin{align*}
\underbar{u}(x) \leq&\ \frac{1}{\left( \frac{\gd}{2} \right)^p} -
\frac{1}{\left( \frac{\gd}{4} \right)^p}\\
=&\ \left( \frac{1}{\gd} \right)^p \left( 2^p - 4^p \right).
\end{align*}
Thus for $p$ chosen large enough with respect to controlled constants one
concludes $\underbar{u} < \phi$ on $U \cap \del B_{\gd}(x_0)$.

It follows that the minimum of $\phi - \underbar{u}$ on $\del U$ is zero and
occurs at $x_0$.  It follows that the normal derivative of $\phi - \underbar{u}$
is nonnegative, and therefore we conclude
\begin{gather*}
\N_i \N_n u \geq - C
\end{gather*}
However, using Lemma \ref{TNlemma1} it is clear that the same argument applies
to $- \phi$, and therefore the result follows.
\end{proof}
\end{lemma}

\begin{lemma} \label{NNC2Estimate} There exists a constant $C$ depending on
$\sup_{0 \leq t \leq 1} \brs{u_t}_{C^1}$,\\ $\sup_{0 \leq t \leq 1} \sup_{x \in
\del M} \brs{\N_i \N_j u_t(x)}$, and $\sup_{0 \leq t \leq 1} \sup_{x \in \del M}
\brs{\N_i \N_n u_t(x)}$ such that for all $x_0 \in \del M$, for all $0 \leq t
\leq 1$ we have
\begin{align*}
\brs{\N_{n} \N_{n} u(x_0)} \leq C.
\end{align*}
\begin{proof} Orthogonally decompose the matrix $W_t$ at $x_0$ in terms of $n$
and $e_i$.  Using the assumed bounds this yields
\begin{gather}
W = \left(
\begin{matrix}
(n-1) u_{nn} & 0\\
0 & u_{nn} g_{|\del M}
\end{matrix}
\right) + \mathcal O(1)
\end{gather}
It is clear that there exists a constant $C > > 0$ such that if $|u_{nn}| > C$
then
\begin{align*}
\det W > \frac{1}{10} C^n > 1 > e^{2n u_t}
\end{align*}
which contradicts the equation $(\star_t)$.  Thus $|u_{nn}(x_0)| < C$ and the
result follows.
\end{proof}
\end{lemma}

\begin{prop} \label{C2Estimate} There exists a constant $C$ such that for all $0
\leq t \leq 1$ we
have
\begin{align*}
\brs{u_t}_{C^2} \leq C
\end{align*}
\begin{proof} By Proposition \ref{C1Estimate} and Lemmas \ref{TTC2Estimate},
\ref{TNC2Estimate}, and \ref{NNC2Estimate} we conclude uniform global $C^1$
bounds and boundary $C^2$ estimates.  Suppose that the maximum of $\brs{\N^2 u}$
occurs at a point in
the interior.  One may follow the calculation of \cite{GV} Proposition 5.1,
which is
justified at any interior point of $M$, to yield the a-priori $C^2$ estimate. 
The result follows.
\end{proof}
\end{prop}

\noindent We can now give the proof of Theorem \ref{dirichlet}.
\begin{proof} As discussed in section 2 it suffices to solve our equation
$(\star_1)$.  Lemma \ref{openness} yields that the set $\Omega$ of $t$ where
$(\star_t)$ is solvable is open in $[0, 1]$.  Proposition \ref{C2Estimate}
yields a uniform $C^2$ estimate for $u_t$.  Thus $(\star_t)$ becomes a uniformly
elliptic equation, and the Evans-Krylov estimates yield uniform
$C^{2,\ga}$ bounds on $u_t$.  Now the Schauder estimates apply to yield uniform
$C^4$ bounds.  Thus $\Omega$ is closed in $[0, 1]$, and hence $(\star_1)$ is
solvable, which completes the proof of existence.  By Proposition \ref{mp}, the
solution to $(\star_1)$ is unique.

To solve the analogous $\sigma_k$ problem, $k < n$, it suffices to observe that
the subsolution we construct for the determinant equation is also a subsolution
to the $\sigma_k$ problem by MacLaurin's inequality.  The proof is otherwise
followed line for line to yield this result.
\end{proof}

\begin{rmk} \label{bndrycond} Note that the proof works equally well if we
require the boundary condition $u_{|\del M} \equiv j$ for any constant $j$. 
Indeed, with minor modification the proof applies to the general boundary value
problem.
\end{rmk}

\section{Proof of Theorem \ref{complete}}
We begin by proving a lemma which serves as a subsolution for solutions to the
Dirichlet problem with large boundary condition.

\begin{lemma} \label{universalsubsoln} Let $(M, \del M, g)$ be a manifold with
boundary.  There exists a function $w$ which is smooth in a controlled
neighborhood of the boundary such that the following holds:  Fix $0 < \ge < <
1$, and let $u$ be a solution to
\begin{align*}
F_1(u) =&\ 0\\
{u}_{|\del M} \equiv&\ - \ln \ge.
\end{align*}
Let $r$ denote distance from the boundary.  Then
\begin{align*}
u \geq - \ln (r + \ge) + \frac{1}{2} \ln(n-1) + w.
\end{align*}
where $w \leq 0$ and $w_{|\del M} = 0$.
\begin{proof} Let $r$ denote distance from $\del M$.  Fix a small
constant $\gd > 0$,
constants $A, p > 0$ and let
\begin{align*}
w =&\ A \left( \frac{1}{(r + \gd)^p} - \frac{1}{\gd^p} \right).
\end{align*}
This choice of $w$ is obviously modeled on our subsolution from section 3,
except that here one thinks of $r + \gd$ as distance from an exterior copy of
$\del M$ instead of distance from a point outside of $\del M$.  The constant
$\gd$ remains fixed for arbitrarily small values of $\ge$.

Fix a small constant $\ge > 0$, and let $\til{r} = r + \ge$.  Again, $\til{r}$
should be thought of as distance from an exterior copy of $\del M$, but this
time one getting arbitrarily close to the actual boundary $\del M$.  Let
\begin{align*}
\underbar{u} =&\ - \ln \til{r} + \frac{1}{2} \ln (n-1) + w.
\end{align*}
Our goal is to show that $\underbar{u}$ is a subsolution for $(\star_t)$ with
boundary condition $u_{|\del M} = - \ln \ge$.  The estimate will proceed in two
steps.  First we estimate $\underbar{u}$ in a small collar neighborhood of the
boundary, where the $-\ln \til{r}$ term dominates the behaviour of
$\underbar{u}$.  Next we exploit the $w$ term using estimates similar to those
in section 3 to control $\underbar{u}$ on the rest of the manifold.

Fix a point $x_0 \in M \backslash \del M$ at which $r$ is smooth, and choose
coordinates at $x_0$ as follows: Let $e_1 =
\frac{\frac{\del}{\del r}}{\brs{\frac{\del}{\del r}}}$, and let $\{e_2, \dots,
e_n \}$ be chosen so that $\{e_i \}$ is an orthonormal basis at $x_0$.  First
observe the preliminary calculation
\begin{align*}
(n-2) \left( \brs{d \underbar{u}}^2 g - d \underbar{u} \otimes d \underbar{u}
\right) =&\ (n-2) \frac{ \left( \til{r} w' - 1 \right)^2}{\til{r}^2} \left(
\begin{matrix}
0 & & & \\
& 1 & & \\
& & \ddots & \\
& & & 1
\end{matrix} \right)
\end{align*}
Therefore at $x_0$ we conclude
\begin{align*}
\hat{\rho} =&\ \rho + (n-2) \N^2 \underbar{u} + \gD \underbar{u} g + (n-2)
\left(\brs{d \underbar{u}}^2 g - d\underbar{u} \otimes d\underbar{u}
\right)\\
=&\ \rho + (n-2) \N^2 w + \gD w g - \frac{1}{\til{r}} \left( (n-2) \N^2 r + \gD
r g
\right)\\
&\ + \frac{1}{\til{r}^2} \left( 
\begin{matrix}
(n-1) & & & \\
& 1 & & \\
& & \ddots & \\
& & & 1
\end{matrix} \right) + (n-2) \frac{ \left( \til{r} w' - 1 \right)^2}{\til{r}^2}
\left(
\begin{matrix}
0 & & & \\
& 1 & & \\
& & \ddots & \\
& & & 1
\end{matrix} \right)\\
=&\ \rho + (n-2) \N^2 w + \gD w g + (n-2) \left( w' \right)^2 \left(
\begin{matrix}
0 & & & \\
& 1 & & \\
& & \ddots & \\
& & & 1
\end{matrix} \right)\\
&\ + \frac{1}{\til{r}^2} \left[ (n-1) g - 2 (n-2) \til{r} w' \left(
\begin{matrix}
0 & & & \\
& 1 & & \\
& & \ddots & \\
& & & 1
\end{matrix} \right) - \til{r} \left( (n-2) \N^2 r + \gD r g \right) \right].
\end{align*}

We now show that the determinant of the bracketed term above, call it $\Phi$, is
positive in a
collar neighborhood of $\del M$ of a fixed width $\eta > 0$.  Observe that when
$\til{r}$
is small this term is dominating the behaviour of $\hat{\rho}$.  We initially
choose $\eta$ small so that the hypersurfaces $\{r = c \}$ are smooth for $c
\leq \eta$.  First note that
$\N^2 r$
is simply the second fundamental form of a smooth hypersurface $\{r = c \}$.
 Therefore it is a tensor with a uniform bound as $r \to 0$ depending only on
$g$.  In particular for $\eta$ chosen small with respect to constants depending
on $g$ we can conclude
\begin{align*}
(n-2) \N^2 r + \gD r g \geq - \gl g.
\end{align*}
Also we can directly compute that on the collar neighborhood of radius $\eta$
\begin{align*}
w' = - \frac{A p}{(r + \gd)^{p+1}} \leq - \frac{A p}{(\eta + \gd)^{p+1}} \leq -
A p
\end{align*}
provided $\eta + \gd < 1$, which is easily arranged.  We therefore conclude that
if we choose $A_0$ and $p$ large with respect to $\gl$, then for any $A \geq
A_0$ and $p \geq p_0$ we have
\begin{align*}
\Phi \geq&\ (n-1) g + \til{r} \left( \begin{matrix}
- \gl & & & \\
& \gl & & \\
& & \ddots & \\
& & & \gl
\end{matrix} \right)
\end{align*}
It follows that if $\eta$ is chosen small with respect to $\gl$, one has $\det
\Phi \geq (n-1)^n$ for $\til{r} \leq 2 \eta$.  Since we have chosen $\ge < < 1$
is follows that for $r < \eta$ we have
\begin{align*}
\det (\hat{\rho}) \geq \frac{(n-1)^n}{\til{r}^{2n}}.
\end{align*}
We would like to show this inequality on the rest of $M$, i.e. for any $r >
\eta$.  Recall from section 3 that given any constant $C > 0$ we may choose our
constants $A$ and $p$ such that
\begin{align*}
(n-2) \N^2 w + \gD w g \geq C g.
\end{align*}
Since $\N^2 r$ is a bounded tensor as described above, we may therefore choose
$A$ and $p$ large so that
\begin{align*}
(n-2) \N^2 w + \gD w g \geq&\ -\rho + \frac{1}{\eta + \ge} \left( (n-2) \N^2 r +
\gD r g \right)\\
\geq&\ - \rho + \frac{1}{2 \eta} \left( (n-2) \N^2 r + \gD r g \right).
\end{align*}
Note here that the choices of $A$ and $p$ depend on $\eta$.  We were careful
above to ensure that the choice of $\eta$ only depended on lower bounds for $A$
and $p$, therefore we are free to choose them still larger, even with respect to
$\eta$.  It follows that at any point where $r$ is smooth we have
\begin{align*}
\det (\hat{\rho}) \geq \left( \frac{1}{\til{r}^2} (n-1) \right)^n =
\frac{(n-1)^n}{\til{r}^{2n}} = e^{2n (- \ln \til{r} + \frac{1}{2} \ln (n-1))}
\geq&\ e^{2n \underbar{u}} 
\end{align*}
where the last inequality follows since $w \leq 0$.  It follows that
$\underbar{u}$ is a subsolution to $(\star_1)$.  It is clear by our construction
of $\underbar{u}$ that the comparison argument of Lemma \ref{subsolnlemma2}
applies to show that any point where $\underbar{u} - u$ achieves it maximum must
be smooth, and hence the argument of Proposition \ref{mp} applies to show that
$u \geq \underbar{u}$.  The lemma follows.
\end{proof}
\end{lemma}

\begin{lemma} \label{universalsupersoln} Let $(M, \del M, g)$ be a manifold with
boundary.  Let $u$ be a solution to
\begin{align*}
F_1(u) =&\ 0
\end{align*}
with any boundary condition.  Then
\begin{align*}
\lim_{x \to \del M} \left[ u(x) + \ln r(x) - \frac{1}{2} \ln(n-1) \right] \leq
0.
\end{align*}
Furthermore, there given any small constant $R > 0$ there exists a constant
$C(R) > 0$ so that given $x_0 \in M$, $B_{R}(x_0) \subset M$ one has
\begin{align*}
u(x_0) \leq C(R).
\end{align*}
\begin{proof} The proof is an adaptation of an argument in \cite{LN} Theorem 4
to the case where the background geometry is not conformally flat.  We observe
that by the Maclaurin inequality it suffices to find a supersolution to the
equation
\begin{align*}
-S := \sigma_1[- \Ric(e^{2 u})] = n e^{2 u}
\end{align*}
to bound the solution to any of the $\sigma_k$ equations from above.  This
equation is given by
\begin{gather} \label{universalsupersoln1}
- \frac{S_g}{n-1} + 2 \gD u + (n-2) \brs{d u}^2 = \frac{n}{n-1} e^{2 u}.
\end{gather}
We proceed to find a local supersolution to this equation.

Take a point $x_0$ in $M$ with distance $d$ from the boundary. 
Consider a geodesic running from the point on the boundary which is closest to
$x_{0}$,
passing through $x_{0}$, and out a small distance $R$ into the manifold to a
point $z_{0}$.  We will fix a small $R$ and a function $f(t)$ based on the
following.  Ensure that we choose both $d$ and $R$
small enough so that given any such point $z_{0}$ as above a distance $R + d$
from the
boundary, then the geodesics inside $B_R(z_0)$ will intersect the
boundary only once, and on this ball $\gD d^{2}(z_{0},\cdot)\geq1.$  Further we
would like to choose $R$ small enough so that there is a solution
to the differential relation on $[0,R^{2}]$
\begin{align*}
(n-2)\left(  f^{\prime}\right)  ^{2}+2f^{\prime \prime}  & \leq0,\\
f^{\prime}  & > \max_{M}|S|+C(g)\\
f(0)  & =0
\end{align*}
where $S$ as above is the scalar curvature on $M$.  In particular, one may
choose
\begin{align*}
f(t) =&\ \sqrt{t + \ge^2} - \ge
\end{align*}
and then for $t$ and $\ge$ chosen small with respect to fixed constants the
required properties are satisfied.

Let $r$ denote the distance function from the point $z_0$, and define a function
$\bar{u}$ on $B_R(z_0)$ by
\begin{align*}
\bar{u} =&\ -\ln(R^{2}-r^{2})+f(R^{2}-r^{2})+\ln2+\ln R+\frac{1}{2}\ln(n-1)
\end{align*}
One directly computes
\begin{align*}
d\bar{u} =&\ \left( \frac{1}{R^{2}-r^{2}}-f^{\prime}\right)  dr^{2}\\
\N^{2}\bar{u} =& \left( \frac{1}{R^{2}-r^{2}}-f^{\prime}\right) 
\N^{2}r^{2}+\left(
\left(\frac{1}{R^{2}-r^{2}}\right)^{2}+f^{\prime\prime}\right)
dr^{2}\otimes dr^{2}\\
\gD \bar{u}=&\ \left(  \frac{1}{R^{2}-r^{2}}-f^{\prime}\right) \gD
r^{2}+\left(  \left(  \frac{1}{R^{2}-r^{2}}\right)  ^{2}+f^{\prime\prime
}\right)  |d r^{2}|^{2}
\end{align*}
It follows that the left hand side of equation (\ref{universalsupersoln1})
becomes
\begin{align*}
& -\frac{S_{0}}{n-1}+2\left(  \frac{1}{R^{2}-r^{2}}-f^{\prime}\right)
\gD r^{2}+2 \left(  \left(  \frac{1}{R^{2}-r^{2}}\right)  ^{2}%
+f^{\prime\prime}\right)  |d r^{2}|^{2}\\
&\ +(n-2)\left(  \frac{1}{R^{2}%
-r^{2}}-f^{\prime}\right)  ^{2}|d r^{2}|^{2}\\
& =\frac{1}{\left(  R^{2}-r^{2}\right)  ^{2}}\left\{
\begin{array}
[c]{c}%
\left[  2\left(  R^{2}-r^{2}\right)  -2 f^{\prime}\left(  R^{2}-r^{2}\right)^2
\right]  (2n+ \tr K)+2\left[  1+\left(  R^{2}-r^{2}\right)  ^{2}f^{\prime\prime
}\right]  4r^{2}\\
+(n-2)(\left[  1-2\left(  R^{2}-r^{2}\right)  f^{\prime}+\left(  R^{2}%
-r^{2}\right)  ^{2}\left(  f^{\prime}\right)  ^{2}\right]  4r^{2}-\frac{S_{0}%
}{n-1}\left(  R^{2}-r^{2}\right)  ^{2}%
\end{array}
\right\}
\end{align*}
where $K:=\N^{2}r^2-2I.$ \ Now using $(n-2)\left(  f^{\prime}\right)
^{2}+2f^{\prime\prime}\leq0,$ we may continue the calculation to bound the above
expression.  In particular
\begin{align*}
& \leq\frac{1}{\left(  R^{2}-r^{2}\right)  ^{2}}\left\{
\begin{array}
[c]{c}%
\left[  2\left(  R^{2}-r^{2}\right)  -2 f^{\prime}\left(  R^{2}-r^{2}\right)^2
\right]  (2n+ \tr K)+8r^{2}\\
+(n-2)(\left[  1-2\left(  R^{2}-r^{2}\right)  f^{\prime}\right]  4r^{2}%
-\frac{S_{0}}{n-1}\left(  R^{2}-r^{2}\right)  ^{2}%
\end{array}
\right\}  \\
& =\frac{1}{\left(  R^{2}-r^{2}\right)  ^{2}}\left\{
\begin{array}
[c]{c}%
4nR^{2}-4nr^{2}+2 \tr K \left(  R^{2}-r^{2}\right)  -2 f^{\prime}\left(  R^{2}%
-r^{2}\right)^2  \gD r^{2}+8r^{2}\\
+4(n-2)r^{2}-8(n-2)\left(  R^{2}-r^{2}\right)  f^{\prime}r^{2}-\frac{S_{0}}%
{n-1}\left(  R^{2}-r^{2}\right)  ^{2}%
\end{array}
\right\}  \\
& =\frac{1}{\left(  R^{2}-r^{2}\right)  ^{2}}\left\{
\begin{array}
[c]{c}%
4nR^{2}+2 \tr K\left(  R^{2}-r^{2}\right)  -2f^{\prime}\left( 
R^{2}-r^{2}\right)^2
\gD r^{2}\\
-8(n-2)\left(  R^{2}-r^{2}\right)  f^{\prime}r^{2}-\frac{S_{0}}{n-1}\left(
R^{2}-r^{2}\right)  ^{2}%
\end{array}
\right\}  \\
& \leq\frac{1}{\left(  R^{2}-r^{2}\right)  ^{2}}\left\{  4nR^{2}%
-(-2 \tr K+\frac{S_{0}}{n-1}\left(  R^{2}-r^{2}\right)  +2 f^{\prime}(R^2 -
r^2))\left(
R^{2}-r^{2}\right)  \right\}
\end{align*}
where in the last line we used that $f' > 0$ and $\gD r^2 \geq 1$.  Applying the
second defining property of $f$ we conclude that
\begin{align*}
- \frac{S_g}{n-1} + 2 \gD u + (n-2) \brs{d u}^2 \leq&\ 4 n R^2 \frac{1}{(R^2 -
r^2)^2}\\
\leq&\ 4 n R^2 \frac{1}{(R^2 - r^2)^2} e^{2 f}\\
=&\ \frac{n}{n-1} e^{2 u}.
\end{align*}
So, given $x_0$ as above and noting that $\bar{u}$ is infinite on $\del
B_{R}(z_0)$, we apply the maximum principle on this ball to conclude that
\begin{align*}
u(x_{0}) &  \leq \bar{u}(x_{0})\\
=&\ -\ln(R^{2}-(R-d)^{2})+\ln2+\ln R+\frac{1}{2}%
\ln(n-1)+f(R^{2}-r^{2})\\
=&\ -\ln(2Rd-d^{2})+\ln2+\ln R+\frac{1}{2}\ln(n-1)+f(2Rd-d^{2})\\
=&\ -\ln d-\ln(2R-d)+\ln2R+\frac{1}{2}\ln(n-1)+f(2Rd-d^{2})\\
=&\ -\ln d+\frac{1}{2}\ln(n-1)-\ln\frac{(2R-d)}{2R}+f(d(2R-d)).
\end{align*}
Taking the limit as $d$ goes to zero yields the first result of the lemma.  To
see the second statement, note that the supersolution $\bar{u}$ can be
constructed as above on any sufficiently small ball in $M$.  Indeed, given $z_0
\in M$ with $B_R(z_0) \subset M$ with $R$ sufficiently small, the estimates
above yield that
\begin{align*}
u(z_0) \leq \bar{u}(z_0) =&\ \ln \frac{(n-1)^{\frac{1}{2}}}{R} + f(R^2) \leq
C(R).
\end{align*}
\end{proof}
\end{lemma}

\noindent We are now ready to give the proof of Theorem \ref{complete}.

\begin{proof} Our proof is similar in nature to \cite{LN} Theorem 4, and indeed
we will exploit an estimate derived there for our purposes.  We reuse the
notation of the previous sections.  The first step is to construct a solution to
the problem
\begin{gather} \label{infdir}
\begin{split}
F_1(u) =&\ 0\\
\lim_{x \to \del M} u(x) =&\ \infty.
\end{split}
\end{gather}
Remark \ref{bndrycond} guarantees the existence of functions $u_j$ solving
\begin{align*}
F_1(u_j) =&\ 0\\
{u_j}_{|\del M} \equiv&\ j.
\end{align*}
We claim that we can extract a subsequence which converges uniformly on compact
sets to a solution to (\ref{infdir}).  First observe that $u_j \geq u_0$ by
Proposition \ref{mp}.  Furthermore, by the second statement of Lemma
\ref{universalsupersoln} we have
that for given $K \subset M \backslash \del M$,
there exists a constant $C = C(K)$ such that $u_j \leq C(K)$ for all $j \geq 0$,
the constant depending on $d(K, \del M)$.
 Therefore we may apply the interior regularity estimates for solutions to
$F_1(u) = 0$ to conclude uniform $C^l$ bounds for any $l$ on any given compact
subset $K \subset M \backslash \del M$.  Interior regularity for such equations
is well established, and one may see for instance \cite{Guan2} Theorems 2.1,
3.1.
 By the Arzela-Ascoli theorem, we conclude that a subsequence $u_{j_n}$
converges uniformly on compact sets to a function $u_{\infty}$.

To show that $u$ is indeed a solution to (\ref{infdir}), first note that by
Proposition \ref{mp} the sequence $\{u_j\}$ is monotonically increasing. 
Applying Lemma \ref{universalsubsoln} we conclude 
\begin{align*}
u_{\infty} \geq&\ \lim_{j \to \infty} u_j\\
\geq&\ \lim_{j \to \infty} \left[- \ln (r + e^{-j}) + \frac{1}{2} \ln (n-1) + w
\right]\\
=&\ - \ln r + \frac{1}{2} \ln(n-1) + w.
\end{align*}
Therefore
\begin{align*}
\lim_{x \to \del M} u_{\infty} \geq&\ \lim_{x \to \del M} - \ln r + \frac{1}{2}
\ln(n-1) + w\\
=&\ \lim_{x \to \del M} - \ln r + \frac{1}{2} \ln(n-1).
\end{align*}
This in fact yields the precise expected asymptotic lower limit for
$u_{\infty}$. Lemma \ref{universalsupersoln} yields the asymptotic upper limit,
implying the precise asymptotic behaviour of $u_{\infty}$ near the boundary.

Finally, we show uniqueness of the solution.  Suppose one had two solutions $u$
and $v$ to equation (\ref{infdir}).  We have already shown that the asymptotic
limits of $u$ and $v$ are the same at $\del M$.  Therefore, let $r$ denote
distance from the boundary of $M$ and consider let $M_{\gd} = \{r \geq \gd \}$
with boundary $\del M_{\gd} = \{r = \gd \}$ for small $\gd > 0$.  Both $u$ and
$v$ are solutions to $F_1(\cdot) = 0$ on $M_{\gd}$ with nearly equal boundary
conditions.  In particular, given $\ge > 0$ we may choose $\gd$ small such that
$u \leq v + \ge$ on $\del M_{\gd}$.  Furthermore, it is clear that $v + \ge$ is
a supersolution to
$F_1(\cdot) = 0$, therefore by Proposition \ref{mp} we conclude that $u \leq v +
\ge$ on $M_{\gd}$.  Taking the limit as $\ge$ goes to $0$ we see that $\gd \to
0$ as well, and so we conclude that $u
\leq v$.  However, the argument is symmetric hence $v \leq u$ and so $u \equiv
v$.
\end{proof}

\section{Conformally Compact Manifolds}

In this section we give an application of Theorem \ref{complete} to the study of
Poincar\'e-Einstein metrics.  The material is 
inspired by the work of Mazzeo-Pacard, and we will refer to \cite{MP} for many
of the details.  We also change notation for this section to that more commonly
used in the study of Poincar\'e-Einstein metrics.  The following definition of
conformally compact metrics contains the basic setup.

\begin{defn}  Let $\bar{X}^{n+1}$ be a compact manifold with boundary $M^n =
\partial
X^{n+1}$.  A Riemannian metric $g_{+}$ 
defined in the interior $X^{n+1}$ is said to be {\em conformally compact} if
there is 
a nonnegative defining function $\rho \in C^{\infty}(\bar{X}^{n+1})$ with 
\begin{align*} 
 \rho &> 0 \ \mbox{ in  } X^{n+1}, \\
\rho &= 0 \ \mbox{ on  } M^n, \\
|\nabla_g \rho| &\neq 0 \ \mbox{ on  } M^n,
\end{align*}
such that $\bar{g} = \rho^2 g_{+}$ defines a Riemannian metric on
$\bar{X}^{n+1}$, and $\bar{g}$ extends at least continuously to $M$. 
The 
manifold $(M^n,\bar{g})$ is called the {\em conformal infinity} of
$(X^{n+1},g_{+})$. 
\end{defn}

Note that one can obtain other defining functions through
multiplication by a positive function; thus the object naturally associated to a
conformally 
compact manifold is not the metric $\bar{g}$ per se (which depends on $\rho$)
but the
conformal class $[\bar{g}]$ of its conformal infinity.  The curvature 
transformation formulas for conformal metrics automatically imply that any
conformally
compact manifold is {\em asymptotically hyperbolic}: i.e., all the sectional
curvatures of 
$g_{+}$ converge to $-1$ at infinity.  

A conformally compact manifold $(X^{n+1},g_{+})$ satisfying the Einstein
condition 
\begin{align} \label{ECon}
 Ric(g_{+}) = -n g_{+}
\end{align}
is called a {\em Poincar\'{e}-Einstein} (P-E) metric.  The 
canonical example of such a metric is $X^{n+1} = B^{n+1} \subset
\mathbb{R}^{n+1}$,
the unit ball, with $g_{+}$ the hyperbolic metric, $\bar{g} = \frac{1}{4} 
(1-|x|^2) g_{+} = ds^2$ the Euclidean metric, and the conformal infinity is the 
round sphere $(\mathbb{S}^n,g = \bar{g}|_{\mathbb{S}^n})$.  Due to its
connection to the 
AdS/CFT correspondence (see \cite{Witten}) there is an extensive literature on
the subject of 
P-E metrics and their physical/geometric properties. 

The question of the existence of a P-E metric with given conformal infinity 
can be interpreted as an asymptotic Dirichlet problem.  In \cite{MP},
Mazzeo-Pacard explored the connection between the 
existence of P-E metrics and the $\sigma_k$-Yamabe problem.  More precisely, 
let $g_{+}$ be a P-E metric; by (\ref{ECon}) the Schouten tensor is given by
\begin{align*}
 A(g_{+}) = -\frac{1}{2}g_{+},
\end{align*}
so that 
\begin{align} \label{skY}
 \sigma_k[ -A(g_{+})] = \frac{1}{2^k} \binom{n + 1}{k} \equiv \beta_{k,n}.
\end{align}
Therefore, a P-E metric is a solution (indeed, the unique solution) of the 
$k$-Yamabe problem in its conformal class, for all $1 \leq k \leq n+1$.  The 
converse is also true: a conformally compact metric $g_{+}$ satisfying
(\ref{skY})
for all $1 \leq k \leq n+1$ is obviously P-E.  The main result of Mazzeo-Pacard
is a more precise statement of this
equivalence: 

\begin{thm}  \label{MPThm} {\em (See \cite{MP}, Theorems 1, 3)}  Let $\Sigma_k$
denote the set of 
 conformally compact metrics on $X^{n+1}$ with Schouten $\sigma_k$-curvature
equal to
$\beta_{k,n}$.  

\vskip.1in \noindent $(i)$ If $g \in \Sigma_k$, then there is a neighborhood
$\mathcal{U}$ of 
$g$ in the space of conformally compact metrics on $X^{n+1}$ such that
$\mathcal{U} \cap \Sigma_k$ 
is an analytic Banach submanifold of $\mathcal{U}$ (with respect to an
appropriate Banach topology).

\vskip.1in \noindent $(ii)$  In addition, 
\begin{align*}
 \mathcal{E} = \bigcap_{k=1}^{n+1} \Sigma_k,
\end{align*}
where $\mathcal{E}$ is the set of Poincar\'{e}-Einstein metrics.  Hence,
$\mathcal{E}$ 
is a finite intersection of locally closed Banach submanifolds, and in
particular is 
always closed in the space of conformally compact metrics on $X^{n+1}$.  
\end{thm}

This equivalence is not just an algebraic curiosity: as Mazzeo-Pacard point
out, 
the linearization of the P-E condition may have a nontrivial finite dimensional
cokernel,
while 
(as we saw in Section 2), the linearization of the Schouten equations do 
not, at least in the negative cone.  On the other hand, it is important to point
out that
when $k \geq 2$, aside from P-E metrics (and perturbations arising from the
above Theorem) there are no 
general existence results for metrics in $\Sigma_k$.  Indeed, for $k \geq 2$,
given a conformal 
infinity $(M^n,[\bar{g}])$ there may be no conformally compact metrics $g_{+}
\in \Sigma_k$ with 
$\rho^2 g_{+}|_{M^n} = \bar{g}$.  

By contrast, it follows from Theorem \ref{infdir} 
that {\em every} conformally compact manifold $(X^{n+1}, g_{+} =
\rho^{-2}\bar{g})$ 
admits a unique conformal metric $h_k = e^{2w_k} \bar{g}$ with
\begin{gather*}
\sigma_k[ - h_k^{-1} Ric(h_k)] = \til{\gb}_{k, n} > 0
\end{gather*}
where we define the constants
\begin{align} \label{betap}
 \tilde{\beta}_{k,n} = \sigma_k( n g ) =  n^k \binom{n + 1}{k},
\end{align}
as the values of $\sigma_k[-Ric]$ for a Poincar\'{e} Einstein metric
normalized as in 
(\ref{ECon}).  Therefore, it is natural to ask whether 
the results of Mazzeo-Pacard have a counterpart for symmetric functions of the 
Ricci tensor.  The answer turns out to be yes.

\begin{thm} Let $(X^{n+1}, g_+)$ be a conformally compact
manifold.  Let $\Theta_k$ denote the set of 
 conformally compact metrics on $X^{n+1}$ with $\sigma_k[-g_+^{-1} \Ric]$ equal
to
$\tilde{\beta}_{k,n}$.  

\vskip.1in \noindent $(i)$ Given a conformally compact metric $g_{+} = \rho^{-2}
\bar{g}$, and $1 \leq 
k \leq n+1$, there is a unique conformally compact metric $h_k = e^{2w_k}
\bar{g} \in \Theta_k$.

\vskip.1in \noindent $(ii)$  Let $\EE$ denote the space of Poincar\'e-Einstein
metrics.  Then 
\begin{align*}
 \mathcal{E} = \bigcap_{k=1}^{n+1} \Theta_k,
\end{align*}
Hence $\mathcal{E}$ is a finite intersection of locally closed Banach
submanifolds, and in
particular is 
always closed in the space of conformally compact metrics on $X^{n+1}$.
\begin{proof} We sketch the details as the proof of a straightforward adaptation
of the proof of Theorem \ref{MPThm}.  That proof relies 
on the structure of the linearized operator for the Schouten tensor equations,
and 
except for some differences of constants, the linearized operator of the
corresponding 
equations for the Ricci tensor is the same.  To be more precise, let
\begin{align*}
\mathcal{H}_k(g_{+},w) =&\ \sigma_k \big[ -Ric(g_{+}) + (n-2)\nabla^2 w + \Delta
w g_{+} \\
&\ \qquad - (n-2) (dw \otimes dw 
-|dw|^2 g_{+}) \big] - \tilde{\beta}_{k,n} e^{2kw}. 
\end{align*}
Thus, if $\mathcal{H}_k[g_{+},w] = 0$ (and $e^{2w}g_{+}$ is conformally compact)
then 
$e^{2w}g_{+} \in \Theta_k$, and conversely.  If $g_{+}$ is P-E, then the
linearization of $\mathcal{H}_k$ with respect to $w$ is given by 
\begin{align*}
 (\mathcal{L}_{Ric})_k \phi = \tilde{c}_{k,n} \Delta \phi - 2k
\tilde{\beta}_{k,n} \phi,
\end{align*}
where 
\begin{align*}
 \tilde{c}_{k,n} = 2(n-1) n^{k-1} \binom{n}{k}.  
\end{align*}
If we linearize the Schouten tensor equations at a P-E metric, the
operator is given by 
\begin{align*}
 (\mathcal{L}_{A})_k \phi = c_{k,n} \Delta  \phi - 2k \beta_{k,n} \phi,
\end{align*}
where 
\begin{align*}
c_{k,n} = 2^{1-k} \binom{n}{k}.
\end{align*}
The essential feature is that, in both cases, there is one positive and one 
negative indicial root of the associated normal operator.  This makes it
possible 
to choose an appropriate weighted space on which $(\mathcal{L}_{A})_k$ and
$(\mathcal{L}_{Ric})_k$ are Fredholm (see Section 2 of \cite{MP}).  After
setting up the right function
spaces and mappings, both statements of Theorem \ref{MPThm} follow from a
version of the implicit
function theorem in, for example, \cite{MS}.
\end{proof}
\end{thm}

Finally, we note that the above characterization of $\mathcal{E}$ can be
weakened considerably, and 
this involves the introduction of a family of $n$ potentially
interesting conformal invariants.
 
\begin{defn} Let $(X^{n+1},g_{+} = \rho^{-2}\bar{g})$, be a conformally compact
manifold.  Fix $1 \leq k \leq n+1$, and let $h_k =
e^{2w_k}\bar{g}$ be the unique conformally compact metric satisfying
\begin{align*} 
\sigma_k [ - h_k^{-1} \Ric(h_k) ] =&\ \til{\gb}_{k,n}.
\end{align*}
Given $1 \leq k \leq n$, let 
\begin{align} \label{Hdef}
H_{k} = w_k - w_{n+1}.
\end{align}
\end{defn}

\begin{prop} \label{Hprop} Let $(X^{n+1}, g_{+} = \rho^{-2} \bar{g})$ be a
conformally compact manifold, and fix $1 \leq k \leq n$.
\begin{enumerate}
\item{$H_{k}$ is a conformal invariant, that is, the definition above does
depend on the choice of conformal background metric $\bar{g}$.} 
\item {$H_{k} \in C^{\infty}(X^{n+1}) \cap C^0(\bar{X}^{n+1})$}
\item{$H_{k} = 0$ on $\partial X^{n+1} = M^n$}
\item{$H_{k} \geq 0$ in $X^{n+1}$.  Moreover, $H(x_0) = 0$
at some point
$x_0 \in X^{n+1}$ if and only if $H_{k} \equiv 0$ and $g_{+}^1 = g_{+}^{n+1}$
is 
a Poincar\'{e}-Einstein metric.}
\end{enumerate}
\begin{proof} (1) If we change $\bar{g}$ to
$e^{2\phi}\bar{g}$ for some 
$\phi \in C^{\infty}(\bar{X}^{n+1})$, then the corresponding solutions to the
$\sigma_k$ problem are respectively
$w_k + 2 \phi$ and $w_n + 2 \phi$.  Therefore, $H =
H(g_{+})$ is uniquely determined 
by a given conformally compact metric, independent of the choice of defining
function.

(2,3) The solutions $w_k$ are always smooth on the interior.  Moreover, by
Lemmas \ref{universalsubsoln} and \ref{universalsupersoln} they have the same
asymptotic limit at $\del M$, hence $H_{|\del M} \equiv 0$.

(4) By MacLaurin's inequality $w_k$ is a supersolution of the $\sigma_{n+1}$
equation, therefore by Proposition \ref{mp} we conclude that $H \geq 0$, and
$H(x_0) = 0$ if and only if $H$ vanishes identically.  If this is the case, the
characterization of equality in the Newton-MacLaurin inequality yields that
$g_+^1 = g_+^{n+1}$ is a Poincar\'e-Einstein metric.  
\end{proof}
\end{prop}

This proposition says that $\mathcal{E}$ is the set of all conformally compact
metrics
$g_{+}$ for which for some $1 \leq k \leq n$, the function $H_{k}
(g_{+})$ vanishes somewhere, and hence everywhere.  From this perspective, the
conformal invariants $H_{k}$ carry the same information for different choices
of $k$.

\section{Remarks on positive curvature}

\begin{rmk} In the notation of the introduction, let $\mathcal S = \{ S_g
> 0 \}$ where here $S_g$ is the Schouten tensor of $g$.  Recall the equation for
the conformal Schouten tensor.
\begin{gather*}
S_{e^{-2u} g} = S_g + \N^2 u + du \otimes du - \frac{1}{2} \brs{d u}^2 g.
\end{gather*}
If we let $u = \ln w$ this is rewritten as
\begin{gather*}
S_{w^{-2} g} = S_g + \frac{1}{w} \N^2 w - \frac{1}{2} \frac{\brs{d w}^2}{w^2} g
\end{gather*}
Now consider an open set $U$ in $\mathbb R^n$ with nonconvex boundary.  Suppose
one had a function $w > 0$ such that $S_{w^{-2} g} > 0$ and $w_{\del U} \equiv
0$.  Since the metric background metric on $U$ is flat, it follows that $\N^2 w
> 0$.  Since $w_{\del U} \equiv 0$, it follows that $\del U$ is a level set of a
convex function, and as such should be convex, which it is not.  Thus we can
never solve for conformal deformation to positive Schouten tensor with
restricted boundary condition.  However, it is not yet clear if we can solve
without the boundary condition, or whether positive Ricci curvature can be
solved for.
\end{rmk}

However, it is easy to show that on a surface with boundary one can always
deform to a metric of positive scalar curvature.
\begin{prop} \label{surfaces} Let $(M^2, \del M, g)$ be a compact Riemannian
surface
with boundary.  There exists $u \in C^{\infty}(M)$ such that $R(e^{-2u} g) > 0$
and $u_{|\del M} \equiv 0$.
\begin{proof} On a surface one has
\begin{align*}
e^{2u} R(e^{-2u} g) = R(g) - 2 \gD u.
\end{align*}
Therefore the problem reduces to solving the Dirichlet problem
\begin{align*}
- 2 \gD u =&\ 1 - R(g)\\
u \equiv&\ 0 \mbox{ on } \del M.
\end{align*}
Solvability of this equation is a known result (\cite{Aubin} Theorem 4.8).
\end{proof}
\end{prop}
\noindent Therefore deformation to positive curvatures remains elusive, and the
nature of the obstructions, if any, are not clear.  To emphasize the issues
here, we formally ask the following question.
\begin{ques} Given $(M^{n}, \del M, g), n \geq 3$, can we conformally deform $g$
to a metric with positive Ricci curvature? Scalar curvature?  Can we do either
while preserving the induced metric on the boundary?
\end{ques}

\bibliographystyle{hamsplain}

\end{document}